\documentclass[11pt]{article}
\usepackage{amsfonts,amsmath,amssymb,graphicx,color}
\usepackage{psfrag}
\usepackage{pdflscape}

\newcommand{\sgn}{\mbox{{sgn}}}

\numberwithin{figure}{section}
\numberwithin{table}{section}

\newcommand{\bequ}{\begin{equation}}     \newcommand{\eequ}{\end{equation}}
\newcommand{\benn}{\begin{equation*}}    \newcommand{\eenn}{\end{equation*}}
\newcommand{\bbma}{\begin{bmatrix}}      \newcommand{\ebma}{\end{bmatrix}}

\newcommand{\R}{\mathbb{R}}

\newcommand{\bsub}{\begin{subequations}}
\newcommand{\esub}{\end{subequations}}

\newtheorem{thm}{Theorem}[section]

\newtheorem{lem}[thm]{Lemma}

\numberwithin{equation}{section}

\newcommand{\comment}[1]{}

\newcommand{\be}{\begin{equation}}
\newcommand{\ee}{\end{equation}}

\newcommand{\bea}{\begin{eqnarray}}
\newcommand{\eea}{\end{eqnarray}}
\newcommand{\beqa}{\begin{eqnarray}}
\newcommand{\eeqa}{\end{eqnarray}}
\newcommand{\beann}{\begin{eqnarray*}}
\newcommand{\eeann}{\end{eqnarray*}}

\newcommand{\brt}{\begin{flushright}}
\newcommand{\ert}{\end{flushright}}
\newcommand{\bc}{\begin{center}}
\newcommand{\ec}{\end{center}}

\newcommand{\bmat}{\left[ \begin{array}}
\newcommand{\emat}{\end{array} \right]}

\newcommand{\beq}{\begin{equation}}
\newcommand{\eeq}{\end{equation}}

\renewcommand{\Re}{{\mathbb R}}          

\newcommand{\fraction}[2]{\textstyle\frac{#1}{#2}}
\def\pref#1{$(\ref{#1})$}
\newcommand{\bproof}{\begin{description} \item[{\it Proof}.] ~ }
\newcommand{\eproof}{\hspace*{\fill}$\Box$\medskip \end{description}}

\newcommand{\defeq}{\stackrel{\rm def}{=}}

\newcommand{\btab}{\begin{tabbing}
\ \ \= thenn \= thenn \= thenn \= thenn \= thenn \= \kill}
\newcommand{\etab}{\end{tabbing}}
\newcommand{\beqas}{\begin{eqnarray*}}
\newcommand{\eeqas}{\end{eqnarray*}}

\newcounter{algo}[section]
\renewcommand{\thealgo}{\thesection.\arabic{algo}}
\newcommand{\algo}[3]{\refstepcounter{algo}
\begin{center}\begin{figure}[h!]
\framebox[\textwidth]{
\parbox{0.95\textwidth} {\vspace{\topsep}
{\bf Algorithm \thealgo : #2}\label{#1}\\
\vspace*{-\topsep} \mbox{ }\\
{#3} \vspace{\topsep} }}
\end{figure}\end{center}}

\newcounter{prog}[section]

\title{An Inexact Successive Quadratic Approximation Method for Convex L-1 Regularized Optimization}
\author{Richard H. Byrd\thanks{Department of Computer Science, University of Colorado,
        Boulder, CO, USA.  This author was supported by National Science Foundation
        grant DMS-1216554 and Department of Energy grant DE-SC0001774.}\and     
       Jorge Nocedal \thanks{Department of Industrial Engineering and Management Sciences, Northwestern University, 
       Evanston, IL, USA.  This author was supported by National Science Foundation grant DMS-0810213, and by Department 
       of Energy grant DE-FG02-87ER25047.} 
       \and
        Figen Oztoprak\thanks{Istanbul Technical University. This author was supported by Department of Energy grant DE-SC0001774,   
        and by a grant from Tubitak.}
      }

\date{\today}
\begin{document}

\maketitle

\begin{abstract}
We study a Newton-like method for the minimization of an objective function $\phi$ that is the sum of a smooth convex function and an $\ell_1$ regularization term. This method, which is sometimes referred to in the literature as a proximal Newton method, computes a step by minimizing a piecewise quadratic model $q_k$ of the objective function $\phi$. In order to make this approach efficient in practice, it is imperative to perform this inner minimization inexactly. In this paper, we give inexactness conditions that guarantee global convergence and that can be used to control the local rate of convergence of the iteration. Our inexactness conditions are based on a semi-smooth function that represents a (continuous) measure of the optimality conditions of the problem, and that embodies the soft-thresholding iteration. We give careful consideration to the algorithm employed for the inner minimization, and report numerical results on two test sets originating in machine learning. \end{abstract}

\section{Introduction}
\label{intro}
\setcounter{equation}{0}

In this paper, we study an inexact Newton-like method for solving optimization problems of the form
\begin{equation}    \label{l1prob}
        \min_{x \in \R^n} \ \phi(x) = f(x) + \mu \|x\|_1 ,
\end{equation}
where $f$ is a smooth convex function and $\mu >0$ is a (fixed) regularization parameter. The method constructs, at every iteration, a piecewise quadratic model of $\phi$ and minimizes this model \emph{inexactly} to obtain a new estimate of the solution. 

The piecewise quadratic model is defined, at an iterate $x_k$, as
\begin{equation}  \label{quadm}
           q_k(x) = f(x_k) + g(x_k)^T (x-x_k) + \fraction{1}{2} (x -x_k)^T H_k (x-x_k) + \mu \| x \|_1 ,
\end{equation}
where $g(x_k) \defeq \nabla f(x_k)$ and $H_k$ denotes the Hessian $\nabla^2 f(x_k)$ or a quasi-Newton approximation to it. After computing an approximate solution $\hat x$ of this model, the algorithm performs a backtracking line search along the direction $d_k=\hat x - x_k$ to ensure decrease in the objective $\phi$. 

We refer to this method as the \emph{successive quadratic approximation method} in analogy to the successive quadratic programming method for nonlinear programming. This method is also known in the literature as a ``proximal Newton method'' \cite{Patr98,sra2011optimization}, but we prefer not to use the term ``proximal'' in this context since the quadratic term in \pref{quadm} is better interpreted as a second-order model rather than as a term that simply restricts the size of the step. 
The paper covers both the cases when the quadratic model $q_k$ is constructed with an exact Hessian or a quasi-Newton approximation. 

The two crucial ingredients in the inexact successive quadratic approximation method are the algorithm used for the minimization of the model $q_k$, and the criterion that controls the degree of inexactness in this minimization. In the first part of the paper, we  propose an inexactness criterion for the minimization of $q_k$ and prove that it guarantees global convergence of the iterates, and that it can be used to control the local rate of convergence. This criterion is based on the optimality conditions for the minimization of \pref{quadm}, expressed in the form of a semi-smooth function that is derived from the soft-thresholding operator. 

The second part of the paper is devoted to the practical implementation of the method. Here, the choice of algorithm for the inner minimization of the model $q_k$ is vital, and we consider two options: {\sc fista} \cite{fista}, which is a first-order method, and an orthant-based method \cite{andrew2007scalable,dss,figi}. The latter is a second-order method where each iteration consists of an orthant-face identification phase, followed by the minimization of a smooth model restricted to that orthant. The subspace minimization can be performed by computing a quasi-Newton step or a Newton-CG step (we explore both options).  A projected bactracking line search is then applied; see section~\ref{obm}.

Some recent work on successive quadratic approximation methods for problem \eqref{l1prob} include: Hsie et al. \cite{hsieh2011sparse}, where \eqref{quadm} is solved using a coordinate descent method, and which focuses on the inverse covariance selection method; \cite{tanscheinberg} which also employs coordinate descent but uses a different working set identification than \cite{hsieh2011sparse}, and makes use of a quasi-Newton model;  and Olsen et al. \cite{olsen2012newton}, where the inner solver is {\sc fista}. None of these papers address convergence for inexact solutions of the subproblem. Recently Lee, Sun and Saunders \cite{lee2012proximal} presented an inexact proximal Newton method that, at first glance, appears to be very close to the method presented here. Their inexactness criterion is, however, different from ours and suffers from a number of drawbacks, as discussed in section~\ref{algo}.  

Inexact methods for solving generalized equations have been studied by Patricksson \cite{patriksson1998cost}, and more recently by Dontchev and Rockafellar \cite{dontchev2013convergence}.  Special cases of the general methods described in those papers result in inexact sequential quadratic approximation algorithms. Patricksson \cite{patriksson1998cost} presents convergence analyses based on two conditions for controlling inexactness.  The first is based on running the subproblem solver for a limited number of steps. The second rule requires that the residual norm be sufficiently small, but it does not cover the inexactness conditions presented in this paper (since the residual is computed differently and their inexactness measure is is different from ours).  The rule suggested in Dontchev and Rockafellar \cite{dontchev2013convergence} is very general, but it too does not cover the condition presented in this paper.
Our rule, and those presented in \cite{dontchev2013convergence,lee2012proximal}, is inspired by the classical inexactness condition proposed by Dembo et al. \cite{DembEiseStei82}, and reduces to it for the smooth unconstrained minimization case (i.e. when $\mu=0$). 

Another line of research that is relevant to this paper is the  global and rate of convergence analysis for inexact proximal-gradient algorithms, which can be seen as special cases of sequential quadratic approximation without acceleration \cite{schmidtLR2011,SalVil12, tappenden2013inexact}.  The inexactness conditions applied in those papers require that the subproblem objective function value  be $\epsilon$-close to the optimal subproblem objective \cite{schmidtLR2011,tappenden2013inexact}, or that the approximate solution be exact with respect to an $\epsilon$-perturbed subdifferential \cite{SalVil12}, for a decreasing sequence $\{\epsilon\}$.

 Our interest in the successive quadratic approximation method is motivated by the fact that it has not received sufficient attention from a practical perspective, where inexact solutions to the inner problem \eqref{quadm} are imperative. Although a number of studies have been devoted to the formulation and analysis of  proximal Newton methods for convex composite optimization problems, as mentioned above,  the viability of the approach in practice has not been fully explored.

This paper is organized in 5 sections. In section~\ref{algo} we outline the algorithm, including the inexactness criteria that govern the solution of the subproblem \pref{quadm}. In sections~\ref{global} and \ref{local}, we analyze the global and local convergence properties of the algorithm. Numerical experiments are reported in section~\ref{numerical}. The paper concludes in section~\ref{finalr} with a summary of our findings, and a list of questions to explore.

\bigskip\noindent{\em Notation.} In the remainder, we let $g(x_k)= \nabla f(x_k)$, and let $\| \cdot \|$ denote any vector norm. We sometimes abbreviate successive quadratic approximation method as ``SQA method'', and note that this algorithm is often referred to in the literature as the ``proximal Newton method''.

\section{The Algorithm}   \label{algo}
\setcounter{equation}{0}
Given an iterate $x_k$, an iteration of the algorithm begins by forming the model \pref{quadm}, where  $\mu >0$ is a given scalar and $H_k \succ 0$ is an approximation to the Hessian $\nabla^2 f(x_k)$. Next, the algorithm computes an \textit{approximate} minimizer $\hat x$ of the subproblem
\begin{equation}  \label{lassop}
    \min_{x \in \R^n} \ q_k(x) = f(x_k) + g(x_k)^T (x-x_k) + \fraction{1}{2} (x -x_k)^T H_k (x-x_k) + \mu \| x \|_1 .
\end{equation}
The point $\hat x$ defines the search direction $d_k= \hat x - x_k$. The algorithm then performs a backtracking line search along the direction $d_k$ that ensures sufficient decrease in the objective $\phi$. The minimization of \eqref{lassop} should be performed by a method that exploits the structure of this problem.

 In order to compute an adequate approximate solution to \pref{quadm}, we need some measure of closeness to optimality. In the case of smooth unconstrained optimization, (i.e. \pref{l1prob} with $\mu = 0$), the norm of the gradient is a standard measure of optimality, and it is common \cite{DembEiseStei82} to require the approximate solution $\hat x$ to satisfy the condition 
\be \label{smoothcond}
           \| g(x_k)+H_k(\hat x -x_k) \|  \leq \eta_k \|g(x_k)\|, \quad\quad  0< \eta_k <1 .
\ee
The term on the left side of \pref{smoothcond} is a measure of optimality for the model $q_k(x)$, in the unconstrained case.  

For problem \eqref{l1prob}, the length of the iterative soft-thresholding (ISTA) step is a natural measure of optimality. The ISTA iteration is given by  
\be \label{ista} 
    x_{\rm ista}= \arg\min_{x}   \ g(x_k)^T(x-x_k) + {\frac{1}{2 \tau}} \|x-x_k\|^2 + \mu \|x\|_1  ,
\ee
where $\tau>0$ is a fixed parameter. It is easy to verify that $\| x_{\rm ista}- x_k \|$ is zero if and only if $x_k$ is a
solution of problem \eqref{l1prob}. We need to express $\| x_{\rm ista}- x_k \|$ in a way that is convenient for our analysis, and for this purpose we note \cite{milzareksemismooth} that  some algebraic manipulations show that $\| x_{\rm ista} - x_k \| = \tau \|F(x_k)\|$, where 
\begin{equation}   \label{Fdef}
        F(x) = g(x) - P_{[-\mu , \mu ]} ( g(x) - x/\tau) .
\end{equation}
Here $P(x)_{[-\mu, \mu]}$ denotes the component-wise projection of $x$ onto the interval $[ -\mu, \mu ]$,  and $\tau$ is a positive scalar. 

One can directly verify that \eqref{Fdef} is a valid optimality measure by noting that $F(x)=0$ is equivalent to 
the standard necessary optimality condition for \pref{l1prob}:
\begin{eqnarray*}
    g_i(x^*) +\mu =0 & \mbox{for   } i \mbox{    s.t.   }   x_i^* >0 , \\
    g_i(x^*) -\mu =0 & \mbox{for   } i \mbox{    s.t.   }   x_i^* <0 , \\
    -\mu \leq g_i(x^*) \leq \mu  & \mbox{for   } i \mbox{    s.t.   }   x_i^* =0  .
\end{eqnarray*}
For the objective $q_k$ of \eqref{lassop},  this function takes the form
\begin{equation}   \label{Fqdef} 
    F_q (x_k; x)= g(x_k)+H_k(x -x_k)-P_{[-\mu,\mu]}( g(x_k)+H_k(x -x_k) -x/\tau).  
\end{equation}
Using the measures \pref{Fdef} and \pref{Fqdef} in a manner similar to \pref{smoothcond},
leads to the condition $\|F_q(x_k; \hat x)\| \leq \eta_k \|F(x_k)\|$. However, depending on the method used to approximately solve
\pref{lassop}, this does not guarantee that $\hat x -x_k$ is a descent direction for $\phi$. To achieve this, we impose the additional condition
that the quadratic model is decreased at $\hat x$.

\bigskip\noindent\textit{Inexactness Conditions.}  A point $\hat x$ is considered an acceptable approximate solution of subproblem \eqref{lassop} if 
\begin{equation}   \label{inx}
          \|F_q(x_k; \hat x)\| \leq \eta_k \|F_q(x_k; x_k)\| \ \mbox{ and } \ \ q_k(\hat x) < q_k(x_k) , \ 
\end{equation}
 for some parameter $0 \leq \eta_k <1$, where $\| \cdot \|$ is any norm.  (Note  that
 $F_q(x_k; x_k)= F(x_k)$, so that the first condition can also be written as
 $ \|F_q(x_k; \hat x)\| \leq \eta_k \|F(x_k)\| $.)
 
 \bigskip
The method is summarized in Algorithm~2.1.

\newpage
 \algo{xalgo}{Inexact Successive Quadratic Approximation (SQA) Method for Problem \pref{l1prob}}
{
	Choose an initial iterate $x_0$. 
	 
	\noindent Select constants  $\theta \in (0,1/2) $ and $0< \tau < 1$ (which is used in definitions \eqref{Fdef}, \eqref{Fqdef}). \\
	
	\noindent
	\textbf{for} $k = 0, \cdots, $ until the optimality conditions of~\eqref{l1prob} are satisfied:	
\smallskip
\begin{itemize}
\item[1.] Compute (or update) the model Hessian $H_k$ and form the piecewise quadratic model \pref{quadm}; 

\item[2.] Compute an inexact solution $\hat x $ of \pref{quadm} satisfying conditions \eqref{inx}. 

\item[3.] Perform a backtracking line search along the direction $d= \hat x -x_k$: starting with $\alpha=1$, find $\alpha \in (0,1]$ such that
\begin{equation}  \label{decrease}
 \phi(x_k) - \phi(x_k+\alpha d) \geq \theta (\ell_k(x_k)- \ell_k(x_k+\alpha d)) ,
\end{equation}
where $\ell$ is the following piecewise linear model of $\phi$ at $x_k$:
\begin{equation}    \label{ell}
         \ell_k(x)  = f(x_k)+g(x_k)^T(x-x_k)+\mu\|x\|_1 .
\end{equation}
\item[4.] Set $x_{k+1}= x_k + \alpha d$. \\

\end{itemize}    
{\bf end(for)}
}

For now, we simply assume that the sequence $\{\eta_k\}$ in \eqref{inx} satisfies $\eta_k \in [0,1)$, but in section~\ref{local} we show that by choosing  $\{\eta_k\}$ and the parameter $\tau$ appropriately, the algorithm achieves a fast rate of convergence.
  One may wonder whether the backtracking line search of Step 3 might hinder sparsity of the iterates. Our numerical experience indicates  that this is not the case because, in our tests,  Algorithm~\ref{xalgo} almost always accepts the unit steplength ($\alpha=1$).
 
 It is worth pointing out that  Lee et al. \cite{lee2012proximal} recently proposed and analyzed an inexactness
criterion that is similar to the first inequality of \pref {inx}.  
The main difference is that they use the subgradient of $q_k$ on the left side of the inequality,
and both norms are scaled by $H_k^{-1}$. They claim similar convergence results to ours, but a worrying consequence
of the lack of continuity of the subgradient of $q_k$ is that their inexactness condition can fail for vectors $x$ arbitrarily close to the exact minimizer of $q_k$. As a result, their criterion is not an appropriate termination test for the inner iteration.
(In addition, their use of the scaling $H_k^{-1}$ precludes setting $ H_k= \nabla^2 f(x_k)$, except for small or highly structured problems.) 

\section{Global Convergence}   \label{global}
\setcounter{equation}{0}

In this section, we show that Algorithm~2.1 is globally convergent under certain assumptions on the function $f$ and the (approximate) Hessians $H_k$. Specifically, we assume that $f$ is a differentiable function with Lipschitz continuous gradient, i.e., there is a constant $M>0$ such that
\begin{equation}   \label{lips}
       \| g(x) - g(y) \| \leq M \| x - y \|,
\end{equation}
for all $x, y$. We denote by $\lambda_{\min}(H_k)$ and $\lambda_{\max}(H_k)$ the smallest and largest eigenvalues of $H_k$, respectively.

\begin{thm} Suppose that $f$ is a smooth  function that is bounded below and that satisfies \pref{lips}. Let $\{ x_k \}$ be the sequence of iterates generated by Algorithm~2.1, and suppose that there exist constants $0< \lambda \leq \Lambda$ such that the sequence $\{ H_k \}$ satisfies
\[
      \lambda_{\min}(H_k) \geq \lambda >0 \quad\mbox{and} \quad \lambda_{\max}(H_k) \leq \Lambda, 
\]
for all $k$.
Then
\begin{equation}    \label{limit}
       \lim_{k \rightarrow \infty} F(x_k) =0 .
\end{equation}
\end{thm}
\textbf{Proof.} We first show that if $\hat x$ is an approximate solution of \pref{quadm} that satisfies the inexactness conditions \pref{inx}, then there is a constant $\gamma >0$ (independent of $k$) such that for  all $k \in \{0,1, \cdots \}$ 
\begin{equation} \label{ldec}
 \ell_k(x_k) - \ell_k(\hat x)   \geq \gamma \|F(x_k)\|^2,
 \end{equation}
 where $\ell_k$ and $F$ are defined in \pref{ell} and \pref{Fdef}.
To see this, note that by \eqref{inx}
\[
0 > q_k(\hat x)-q_k(x_k) = \ell_k(\hat x)-\ell_k(x_k) + \fraction{1}{2}(\hat x-x_k)^TH_k(\hat x-x_k) ,
\]
and therefore
\begin{equation}    \label{lpre}
      \ell_k(x_k) - \ell_k(\hat x) > \fraction{1}{2}(\hat x-x_k)^TH_k(\hat x-x_k) \geq  \fraction{1}{2}\lambda \|\hat x-x_k\|^2 .
\end{equation}
 
Next, since $F(x_k)= F_q(x_k;x_k)$, and using \pref{inx}  and the contraction property of the projection, we have that
\begin{align*}
  (1-\eta_k)\|F(x_k)\| &= (1-\eta_k)\|F_q(x_k;x_k)\| \\
     & \leq \|F_q(x_k;x_k)\| - \|F_q(x_k,\hat x)\|  \\
     & \leq  \|F_q( x_k;\hat x)-F_q(x_k;x_k)\| \\
  & = \|H_k(\hat x-x_k)-P_{[-\mu,\mu]}(g(x_k) + H_k(\hat x-x_k)-\fraction{1}{\tau}\hat x) + P_{[-\mu,\mu]}(g(x_k)-\fraction{1}{\tau}x_k)\|\\ 
  & \leq \|H_k(\hat x-x_k)\| + \|\fraction{1}{\tau}(\hat x-x_k)-H_k(\hat x-x_k)\| \\
  & \leq \fraction{1}{\tau}\|\hat x-x_k\| + 2\|H_k\|\|\hat x - x_k\| \\
  & = \left(\fraction{1}{\tau} + 2\|H_k\| \right)\|\hat x-x_k\|   \\
  & \leq \left(\fraction{1}{\tau} + 2\Lambda \right)\|\hat x-x_k\|.
\end{align*}
Combining this expression with \pref{lpre}, we obtain \pref{ldec} for
\[
 \gamma = \frac{\lambda}{2}\left( \frac{1-\eta}{\frac{1}{\tau}+2\Lambda}\right)^2.
\]
Note that $\gamma>0$ as $\tau,\lambda, \Lambda>0$ and $\eta\in[0,1)$.

\smallskip
 
Let us define the search direction as $d= \hat x- x_k$. We now show that by performing a line search along $d$ we can ensure that the algorithm provides sufficient decrease in the objective function $\phi$, and this will allow us to establish the limit \eqref{limit}. 

Since $g(x)$ satisfies the Lipschitz condition \pref{lips}, we have
\[
f(x_k+\alpha d) + \mu \|x_k+\alpha d\|_1 \leq f(x_k) + \alpha g(x_k)^Td + \frac{M}{2}\alpha^2\|d\|^2 + \mu\|x_k+\alpha d\|_1 ,
\]
and thus
\[
\ f(x_k) + \mu \|x_k\|_1 - f(x_k+\alpha d) - \mu \|x_k+\alpha d\|_1 \geq - \alpha g(x_k)^Td - \frac{M}{2}\alpha^2\|d\|^2 - \mu\|x_k+\alpha d\|_1 + \mu\|x_k\|_1.
\]
Recalling the definition of $\ell_k$, we have
\[
\phi(x_k) - \phi(x_k+\alpha d) \geq \ell_k(x_k) - \ell_k(x_k+\alpha d) - \frac{M}{2}\alpha^2\|d\|^2.
\]
By convexity of the $\ell_1$-norm, we have that
\[
\ell(x_k) - \ell(x_k+\alpha d) \geq  \alpha(\ell(x_k) - \ell(x_k+d)).
\]
Combining this inequality with \eqref{lpre}, and recalling that $x+d = \hat x$, we obtain for $ \theta\in(0,1)$,
\begin{align}
 \phi(x_k) - \phi(x_k+\alpha d) - \theta(\ell(x_k) - \ell(x_k+\alpha d)) & \geq (1-\theta)(\ell(x_k) - \ell(x_k+\alpha d))- \frac{M}{2}\alpha^2\|d\|^2  \nonumber \\
 & \geq (1-\theta)\alpha(\ell(x_k) - \ell(x_k+d))- \frac{M}{2}\alpha^2\|d\|^2 \nonumber  \\
 & \geq (1-\theta)\alpha\frac{\lambda}{2}\|d\|^2- \frac{M}{2}\alpha^2\|d\|^2\nonumber\\
 & = \fraction{1}{2}\alpha\|d\|^2 \left((1-\theta)\lambda-M\alpha\right) \nonumber \\
 & \geq  0 , \label{spray}
\end{align}
provided $\left(  (1-\theta)\lambda-M\alpha \right) \geq 0$. Therefore, the sufficient decrease condition \eqref{decrease} is satisfied for any steplength $\alpha$ satisfying
\[
0 \leq \alpha \leq (1-\theta)\frac{\lambda}{M},
\]
and if the backtracking line search cuts the steplength in half (say) after each trial, we have that the steplength chosen by the line search satisfies
\[
\alpha \geq  (1-\theta)\frac{\lambda}{2 M} .
\]
Thus, from \eqref{spray} and \eqref{ldec} we obtain 
\[
 \phi(x_k) - \phi(x_k+\alpha d) \geq \theta(1-\theta)\frac{\lambda}{2 M}\gamma\|F(x_k)\|^2.
\]
Since $f$ is assumed to be bounded below, so is the objective function $\phi$, and given that the decrease in $\phi$ is proportional to $\| F(x_k) \|$ we obtain the limit \pref{limit}. 
\hspace*{\fill}$\Box$\medskip\medskip

We note that to establish this convergence result it was not necessary to assume convexity of $f$.

\section{Local Convergence}   \label{local}
\setcounter{equation}{0}

To analyze the local convergence rate of the successive quadratic approximation method, we use
the theory developed in Chapter 7 of Facchinei and Pang \cite{facchinei2003finite}.
To do this, we first show that, if $x^*$ is a nonsingular minimizer of $\phi$, then the functions
$F_q (x; \cdot): \R^n \rightarrow \R^n$
are a family of uniformly Lipschitzian  nonsingular homeomorphisms for all $x$ in a
neighborhood of $x^*$.

\begin{lem}  \label{strongmono}
If $H$ is a symmetric positive definite matrix with smallest eigenvalue $\lambda >0$,
then the function of $y$ given by 
\[ 
           F_q (x;y)= g(x)+H(y-x)-P_{[-\mu,\mu]}( g(x)+H(y-x) -y/\tau),  
\]
is strongly monotone if $\tau < 1/\|H\|$. Specifically, for any vectors $y,z \in \R^n$,
\be \label{smono}
           (z-y)^T(F_q (x;z)- F_q (x;y)) \geq \fraction{1}{2}{\lambda} \|z-y\|^2 .
\ee
\end{lem}
\textbf{Proof.}
%
%
It is straightforward to show that for any scalars $a \neq b$ and interval $[-\mu,\mu]$,
\begin{equation} \label{projlip}  
       0 \leq { P_{[-\mu,\mu]}(a) -  P_{[-\mu,\mu]}(b)  \over a-b }  \leq 1.
\end{equation}
Therefore  for any vectors $y$ and $z$, and for any index $i \in \{1, \cdots, n\}$ we have
\begin{align*}
 F_q (x;z)_i- F_q (x;y)_i &= H(z-y)_i-P_{[-\mu,\mu]}( g_i(x)+H(z-x)_i -z_i/\tau)] \\
                                      & \quad + P_{[-\mu,\mu]}( g_i(x)+H(y-x)_i -y_i/\tau)] \\
                                   &= H(z-y)_i -\bar d_i [H(z-y)_i -(z_i-y_i)/ \tau] ,
\end{align*}
where $\bar d_i \in [0,1]$ is a scalar implied by (\ref{projlip}). This implies that
\[
             F_q (x;z)- F_q (x;y) =H(z-y) + D(\fraction{1}{ \tau}I-H)(z-y) ,
\]
where $D = {\rm diag}(\bar d_i)$. Hence
\begin{equation}    \label{ezza}
(z-y)^T (F_q (x;z)- F_q (x;y))=(z-y)^T H(z-y) + (z-y)^TD(\fraction{1}{ \tau}I-H)(z-y).
\end{equation}
Since the right hand side is a quadratic form, we symmetrize the matrix, and if we let $w=z-y$, the right side is
\begin{equation}  \label{loud}
w^T [H+ \tau^{-1}D - \fraction{1}{2}(DH+HD)]w .
\end{equation}
To show that the symmetric matrix inside the square brackets is positive definite, we note that since
$(\tau H -D)^T(\tau H -D) = \tau^2 H^2 -\tau (HD+DH) +D^2$ is positive semi-definite,
we have that 
\[
w^T(HD+DH)w \leq w^T(\tau H^2 +\tau^{-1}D^2)w.
\] 
Substituting this into \eqref{loud} yields
\begin{align*}
    w^T \left[H+ \tau^{-1} D -\fraction{1}{2}(DH+HD)\right]w & \geq w^T\left[H-{\tau \over 2}H^2 + {D-D^2/2 \over \tau }\right]w \\
          &    \geq w^T\left[H-{\tau \over 2}H^2 \right]w ,
\end{align*}
since  $D-D^2/2$ is positive semi-definite given that the elements of the diagonal matrix $D$ are in $[0,1]$.
If $\lambda_i$ is an eigenvalue of $H$, the corresponding eigenvalue of the matrix $H-{\tau \over 2}H^2$
is $\lambda_i-\tau \lambda_i^2 /2 \geq \lambda_i/2$  since our assumption on $\tau$ implies $1>\tau \|H\|\geq \tau \lambda_i$.
Therefore, we have from \eqref{ezza} that 
\[
(z-y)^T (F_q (x;z)- F_q (x;y)) \geq \fraction{1}{2}\lambda\|z-y\|^2 .
\]
\hspace*{\fill}$\Box$\medskip
\comment{
We claim that if $\tau \|H\| <1$ the matrix $H+ D( {1 \over \tau}I -H)$ is positive definite. 
To show this, let $w=y-z$. For an arbitrary scalar $\delta \in (0,1)$, 
decompose $D=D_1 +D_2$, where the diagonal matrix $D_1$ contains the values $d_i$ in $[0, \delta]$
and $D_2$ contains those in $(\delta, 1]$. Similarly decompose $w=w_1+w_2$.
Then
\begin{eqnarray*}
 w^T Hw +w^TD({1\over \tau}I-H)w &=&  w^T(I-D_1)Hw + w^T(\fraction{1}{\tau}D -D_2 H)w \\
 &\geq& (\lambda -\delta \|H\|)\|w_1\|^2 + w_2^TD_2({1 \over \tau}-\|H\| ) w_2 \\
 &\geq& (\lambda -\delta \|H\|) \|w_1\|^2 + ({1\over \tau}-\|H\| )\delta \|w_2\|^2 \\
 &\geq&  \min\{ \lambda -\delta \|H\| , ({1 \over \tau}-\|H\|) \delta \} \|w\|^2.
\end{eqnarray*}
If the value $\tau < 1/ \|H\|$, and we choose $\delta = \tau \lambda$
then $F_q$ is strongly monotone with constant 
\[ \lambda(1 -\tau \|H\| ). \]
If $\tau < {1 \over 2 \|H\|}$
the constant is $\lambda/2$.
\hspace*{\fill}$\Box$\medskip
}   

Inequality \eqref{smono} establishes that $F_q(x;\cdot)$ is strongly monotone. 
Next we show that, when $H$ is defined as the Hessian of $f$,  the functions $ F_q (x;\cdot)$ are homeomorphisms and that they represent an accurate
approximation to the function $F$ defined in \eqref{Fdef}.

\begin{thm} \label{homeo} If $\nabla^2 f(x^*)$ is positive definite and $\tau < 1/\|\nabla^2 f(x^*) \|$, then
there is a neighborhood $\cal{N}$ of $x^*$ such that 
for all $x \in \cal{N}$ the functions of $y$ given by
\begin{equation}   \label{ella} 
      F_q (x;y)= g(x) +\nabla^2 f(x)(y-x)-P_{[-\mu,\mu]}( g(x)+ \nabla^2 f(x)(y-x) -y/\tau)
 \end{equation}
are a family of homeomorphisms from $\R^n$ to $\R^n$, whose inverses
$F_q^{-1}(x;\cdot) $ are uniformly Lipschitz continuous.
In addition, if $\nabla^2 f(x)$ is Lipschitz continuous, then 
there exists a constant $\beta>0$ such that 
\be \label{fqerror}
\| F(y) -F_q(x;y) \| \leq \beta \|x-y\|^2  
\ee
for any $x \in \cal{N}$.
\end{thm}
\textbf{Proof.}
Since $\nabla^2 f(x)$ is continuous, there is a neighborhood ${\cal N}$ of $x^*$ and a positive constant
$\lambda$ such that
$\lambda_{min} (\nabla^2 f(x)) \geq \lambda >0$ and $\tau \|\nabla^2 f(x) \| < 1$,
for all $x\in \cal{N}$.  
It follows from Lemma \ref{strongmono} that for any such $x$,
the function $F_q(x;y)$ given by \eqref{ella} is strongly (or uniformly) monotone with constant greater than $ \lambda/2$. We now invoke the Uniform Monotonicity Theorem (see e.g.  Theorem~6.4.4 in \cite{OrteRhei70}), which states that if a function $F:\R^n \rightarrow \R^n$ is continuous and uniformly monotone, then $F$ is a homeomorphism of $\R^n$ onto itself.
We therefore conclude that  $F_q(x;y)$ is a homeomorphism.  

In addition, we have from \pref{smono} and the Cauchy-Schwartz inequality that
\[ 
     \|z-y\| \| F_q (x;z)- F_q (x;y)\| \geq (z-y)^T(F_q (x;z)- F_q (x;y)) \geq  \fraction{1}{2} \lambda \|z-y\|^2 ,
\]
which implies Lipschitz continuity of $F_q^{-1}(x;\cdot) $ with constant $ 2 / \lambda $.
To establish \pref{fqerror}, note that
\begin{align*}
F(y)-F_q(x;y) & = g(y) -P_{[-\mu,\mu]}( g(y) -y/\tau) - \left(g(x) +\nabla^2 f(x)(y-x) \right)  \\
   & +P_{[-\mu,\mu]}( g(x)+\nabla^2 f(x)(y-x) -y/\tau),
\end{align*}
and thus
\begin{eqnarray*}
\|F(y)-F_q(x;y)\| &\leq& \| g(y) - (g(x) +\nabla^2 f(x)(y-x)) \| \\
         & & +
          \|P_{[-\mu,\mu]}( g(y) -y/\tau) - P_{[-\mu,\mu]}( g(x)+\nabla^2 f(x)(y-x) -y/\tau) \| \\
       &\leq& 2 \| g(y) - ( g(x) +\nabla^2 f(x)(y-x)) \|   \\
       & = & O(\|y-x\|^2) ,
\end{eqnarray*}
by the non-expansiveness of a projection onto a convex set and Taylor's theorem. 
\hspace*{\fill}$\Box$\medskip

Theorem~\ref{homeo} shows that $F_q(x;y)$ defines a strong nonsingular Newton approximation 
in the sense of Definition~7.2.2 of Pang and Facchinei \cite{facchinei2003finite}. This implies quadratic convergence for the (exact) successive quadratic approximation (SQA) method.

\begin{thm} If $\nabla^2 f(x)$ is Lipschitz continuous and positive definite at $x^*$, and $\tau < 1/\|\nabla^2 f(x^*) \|$, then
there is a neighborhood of $x^*$ such that, if $x_0$ lies in that neighborhood, the iteration that
defines $x_{k+1}$ as the unique solution to
\[ F_q (x_k;x_{k+1})= 0 \]
converges quadratically to  $x^*$.
\end{thm}
\textbf{Proof.}
By Theorem \ref{homeo}, $F_q(x_k;y)$ satisfies the definition of a nonsingular strong Newton
approximation of $F$ at $x^*$, given by Facchinei and Pang (\cite{facchinei2003finite}, 7.2.2) and thus by Theorem~7.2.5
of that book the local convergence is quadratic.
\hspace*{\fill}$\Box$\medskip

Now we consider the inexact SQA algorithm that, at each step, computes a point $y$
satisfying
\be
 F_q (x_k;y)= r_k ,
\ee
where $r_k$ is a vector such that $\|r_k\| \leq \eta_k \|F(x_k)\|$ with $\eta_k <1$; see \eqref{inx}.
We obtain the following result for a method that sets $x_{k+1}=y$. 

\begin{thm} \label{inexact} Suppose that $\nabla^2 f(x)$ is Lipschitz continuous and positive definite at $x^*$, $\tau < 1/\|\nabla^2 f(x^*) \|$,
and that $x_{k+1}$ is computed by solving
\[ F_q (x_k;x_{k+1})= r_k , \]
where $\|r_k\| \leq \eta_k \|F(x_k)\|$.
Then, there is a neighborhood $\cal{N}$ of $x^*$ and a value $\bar{\eta}>0 $ such that if $\eta_k \leq \bar{\eta} $
for all $k$ and if $x_0 \in \cal{N}$ then the sequence $\{x_k\}$ converges Q-linearly to $x^*$.
In addition if $\eta_k \rightarrow 0$, then the convergence rate of $\{x_k\}$ is Q-superlinear.
Finally, if for some $\tilde{\eta}$, $\eta_k \leq \tilde{\eta} \|F(x_k)\|$ then the convergence rate
is Q-quadratic.  
\end{thm}
\textbf{Proof.}
By Theorem \ref{homeo}, the iteration described in the statement of the theorem satisfies
all the conditions of Theorem~7.2.8 of \cite{facchinei2003finite}. The results then follow immediately from that
theorem.
\hspace*{\fill}$\Box$\medskip

We have shown above the the inexact successive quadratic approximation (SQA) method with $\alpha_k =1$ yields a fast rate of convergence. We now show that this inexact SQA algorithm will select the steplength $\alpha_k=1$ in a neighborhood of the solution.  In order to do so, we strengthen the inexactness conditions \eqref{inx} slightly so that they read
\begin{equation}   \label{inxx}
          \|F_q(x_k; \hat x)\| \leq \eta_k \|F_q(x_k; x_k)\| \ \mbox{ and } \ \ q_k(\hat x)-q_k(x_k) \leq 
          \zeta (\ell_k(\hat x)-\ell_k(x_k) )  ,
          \end{equation}
where $\eta_k <1$,  $\zeta \in (\theta , 1/2)$ and $\theta$ is the input parameter of Algorithm~2.1 used in \eqref{decrease}. Thus, instead of simple decrease, we now impose sufficient decrease in $q_k$. 

\begin{lem} \label{nloa}
 If $H_k$ is positive definite, the inexactness condition \pref{inxx} is satisfied by any sufficiently accurate solution to \pref{lassop}.
\end{lem}
\textbf{Proof.}
If we denote by $\bar{y}$ the (exact) minimizer of $q_k$, we claim that 
\begin{equation}  \label{masia}
    q_k(x_k) - q_k(\bar{y}) \geq \fraction{1}{2} (\ell_k(x_k) - \ell_k(\bar{y}) ). 
\end{equation}
To see this, note that $q_k(y) = \ell_k(y) + {1\over2} (y-x_k)^TH_k (y-x_k)$, 
and since $\ell_k$ and $q_k$ are convex and $\bar{y}$ minimizes $q_k$, there exists a
vector $v \in \partial \ell_k(\bar{y})$ such that $ v+ H_k (\bar{y}-x_k) =0$.
By convexity
\begin{equation}  \label{ficus}
   \ell_k(x_k) \geq \ell_k(\bar{y}) + v^T(x_k - \bar{y} )= \ell_k(\bar{y}) + (\bar{y}-x_k)^TH_k(\bar{y}-x_k).
\end{equation}   
Therefore, 
\[ q_k(x_k) - q_k(\bar{y}) = \ell_k(x_k)- \ell_k(\bar{y}) - \fraction{1}{2}(\bar{y}-x_k)^TH_k (\bar{y}-x_k)  \\ 
                    \geq \fraction{1}{2} (\ell_k(x_k)- \ell_k(\bar{y}) ) ,
\]
which proves \eqref{masia}. 

Now consider the continuous function $q_k(x)- q_k (x_k) - \zeta (\ell_k(x) - \ell_k(x_k) )$.
By \pref{masia} its value at $x=\bar{y}$ is 
\be
   q_k(\bar y )- q_k (x_k) - \zeta ( \ell_k(\bar y ) -  \ell_k(x_k) )    \leq (\fraction{1}{2} -\zeta)( \ell_k(\bar y ) -  \ell_k(x_k) ) <0,
\ee
where the last inequality follows from \eqref{ficus}.
Therefore by continuity, the value of this function for any $x$ in some neighborhood of $\bar{y}$ is negative, implying that \pref{inxx} is satisfied by
any approximate solution $\hat{x}$ sufficiently close to $\bar{y}$.
\hspace*{\fill}$\Box$\medskip

\begin{thm} \label{approxaccept}
Suppose that $H_k=\nabla^2f(x_k)$ in Algorithm~2.1, and that we modify Step 2 in that algorithm to require that the approximate solution $\hat x$ satisfies \eqref{inxx} instead of \eqref{inx}. If we assume that $\nabla^2f(x)$ is Lipschitz continuous, then for all $k$ sufficiently large we have $\alpha_k =1$.
\end{thm}
\textbf{Proof.}
Given that $\hat x = x_k +d_k$ satisfies \pref{inxx}, it follows from
Taylor's theorem, the Lipschitz continuity of $\nabla^2f(x)$, and equation \pref{lpre} that for some constant $\rho>0$
\begin{eqnarray*}
\phi(x_k + d_k) -\phi(x_k) &=&  [ \phi(x_k + d_k) -\phi(x_k) - q_k(x_k+d_k) +q_k(x_k)] \\
               &   & -(q_k(x_k)-q_k(x_k+d_k)) \\
               &\leq&  -\zeta (\ell_k(x_k)- \ell_k(x_k +d_k) ) +\rho \|d_k\|^3     \\
               &\leq&  \theta ( \ell_k(x_k + d_k)- \ell_k(x_k) ) + (\zeta - \theta) ( \ell_k(x_k+ d_k)- \ell_k(x_k) ) + \rho\|d_k\|^3     \\
               &\leq&  \theta ( \ell_k(x_k+ d_k)- \ell_k(x_k) ) - (\zeta - \theta) {\lambda \over 2} \|d_k\|^2 + \rho \|d_k\|^3     \\
               &\leq&  \theta ( \ell_k(x_k + d_k)- \ell_k(x_k) )
\end{eqnarray*}
if $\|d_k\| \leq (\zeta -\theta) \lambda/2\rho$ .
Since the global convergence analysis implies $\|d_k\| \rightarrow 0$, we have from \eqref{decrease} that eventually the steplength $\alpha_k =1$ is accepted and used.
\hspace*{\fill}$\Box$\medskip

\bigskip
We note that \eqref{inxx} is stronger than \eqref{inx}, and therefore, all the results presented in this and the previous section apply also to the strengthened condition \eqref{inxx}.
Theorem \ref{approxaccept} implies that if Algorithm~2.1 is run with the strengthened accuracy condition \pref{inxx}, 
and $H_k=\nabla^2f(x_k)$, then once the iterates are close enough to 
a nonsingular minimizer $x^*$, the iterates have the linear, superlinear or quadratic convergence rates described in Theorem~\ref{inexact} if
$\eta_k$ is chosen appropriately.

\section{Numerical Results}   \label{numerical}
\setcounter{equation}{0}
One of the goals of this paper is to investigate whether the successive quadratic approximation ({\sc sqa}) method is, in fact, an effective approach for solving convex $\ell_1$ regularized problems of the form \eqref{l1prob}. Indeed, it is reasonable to ask whether it might be more effective to apply an algorithm such as {\sc ista} or {\sc fista}, directly to problem \eqref{l1prob}, rather than performing an inner iteration on the subproblem \eqref{lassop}.  Note that each iteration of {\sc fista} requires an evaluation of the gradient of the objective \eqref{l1prob}, whereas each inner iteration for the subproblem \eqref{lassop} involves the product of $H_k$ times a vector.

  To study this question, we explore various algorithmic options within the successive quadratic approximation method, and evaluate their performance using data sets with different characteristics. One of the data sets concerns the covariance selection problem (where the unknown is a matrix), and the other involves a logistic objective function (where the unknown is a vector). Our benchmark is {\sc fista} applied directly to problem \eqref{l1prob}. {\sc fista} enjoys convergence guarantees when applied to problem \eqref{l1prob}, and is generally regarded as an effective method. 

We employ two types of methods for solving the subproblem \eqref{lassop} in the successive quadratic approximation method: {\sc fista} and an orthant based method {\sc (obm)}  \cite{andrew2007scalable,dss,figi}. The orthant based method (described in detail in section~\ref{obm}) is a two-phase method in which an \emph{active orthant face} of $\R^n$ is first identified, and a subspace minimization is then performed with respect to the variables that define the orthant face. 
The subspace phase can be performed by means of a Newton-CG iteration, or by computing a quasi-Newton step; we consider both options. 

\bigskip\noindent
The methods employed in our numerical tests are as follows.

\begin{itemize}
\item[ ] {\sc \textbf{FISTA.}} This is the {\sc fista} algorithm \cite{fista} applied to the original problem \eqref{l1prob}. We used the
                               implementation from the  
                               TFOCS package,  called N83 \cite{becker2011templates}.  This implementation differs from the 
                               (adaptive) algorithm described by 
                                Beck and Teboulle~\cite{fista} in the way the Lipschitz parameter is updated, and
                                performed significantly better in our test set than the method in \cite{fista}. 
\item[ ] \textbf{PNOPT.}  This is the sequential quadratic approximation (proximal Newton) method of Lee, Sun and Saunders \cite{lee2012proximal}.  The Hessian $H_k$ in the subproblem \eqref{lassop} is updated using the limited memory BFGS formula, with a (default) memory of 50. (The {\sc pnopt}
package also allows for the use of the exact Hessian, but since this matrix must be formed and factored at each iteration, its use is impractical.) The subproblem \eqref{lassop} is solved using the N83 implementation of {\sc fista} mentioned 
                             above. {\sc pnopt} provides the option
                            of using {\sc sparsa} \cite{sparsa} instead of N83 as an inner solver, but the performance of {\sc sparsa} was not robust in our tests, and we will 
                            not report results with it.
\item[ ] \textbf{SQA.}  Is the sequential quadratic approximation method described in Algorithm~\ref{xalgo}. We implemented 3 variants that differ in the method used to solve the subproblem \eqref{lassop}.
   \begin{itemize}
       \item[ ] \textbf{SQA-FISTA.} This is an {\sc sqa} method using {\sc fista-n}{\small 83} to solve the subproblem \eqref{lassop}. The 
                  matrix $H_k$ is the exact Hessian $\nabla^2 f(x_k)$; each inner {\sc fista} iteration requires two multiplications with $H_k$.
               \item[ ] \textbf{SQA-OBM-CG.}  This is an {\sc sqa} method that employs an orthant based method 
                       to solve the subproblem \eqref{lassop}. The {\sc obm} method performs the subspace minimization step using  a Newton-CG iteration. The number of CG iterations varies
                       during the course of the (outer) iteration according to the rule $\min\{3, 1+\lfloor k/10\rfloor \}$, where $k$ is the outer iteration number.
              \item[ ] \textbf{SQA-OBM-QN.}  This is an {\sc sqa} method where the inner solver is an {\sc obm} method in which the
                      subspace phase consists of a limited memory BFGS step, with a memory of 50.  The correction
                      pairs used to update the quasi-Newton matrix employ gradient differences from the outer iteration (as in
                      {\sc pnopt}).
   \end{itemize}
\end{itemize}

The initial point was set to the zero vector in all experiments, and the iteration was terminated if $\| F (x_k)\|_\infty \leq 10^{-5}$, where $F$ is defined in \eqref{Fdef}. The maximum number of outer iterations for all solvers was $3000$. In the SQA method, the parameter  $\eta_k$ in  the inexactness condition \eqref{inx} was defined as $\eta_k = \max \{1/k , 0.1\}$, and we set $\theta=0.1$ in \eqref{decrease}. For {\sc pnopt} we set  `ftol'=$1e-16$, and `xtol'=$1e-16$ (so that those two tests do not terminate the iteration prematurely), and chose `Lbfgs\_mem'=50. 

We noted above that Algorithm~\ref{xalgo} can employ the inexactness conditions \eqref{inx} or \eqref{inxx}. We implemented both conditions, with  $\zeta=\theta=0.1$, and obtained identical results in all our runs. 

We now describe the numerical tests performed with these methods.

\subsection{Inverse Covariance Estimation Problems } 
The task of estimating a high dimensional sparse inverse covariance matrix is closely tied to the topic of Gaussian Markov random fields \cite{picka:06}, and arises in a variety of recognition tasks. This model can be used to recover a sparse social or genetic network from user or experimental data.

A popular approach to solving this problem \cite{Banerjee:08,Banerjee:06} is to minimize the negative log likelihood function, under the assumption of normality, with an additional $\ell_1$ term to enforce sparsity in the  estimated inverse covariance matrix. We can write the optimization problem as
\begin{equation}    \label{cov}
        \min_{P \in \R^{n\times n} } \ {\rm tr}(S P) - \log \det P +\mu \| P\| ,
\end{equation}
where $S$ is a given sample covariance matrix, $P$ denotes the unknown inverse covariance matrix, $\mu$ is the regularization parameter, and $\| P \| \defeq \| vec( P) \|_1$.  
We note that the Hessian of the first two terms in \eqref{cov} has a very special structure: it is given by $P^{-1} \otimes P^{-1}$.

Since the objective is not defined when $\det(P)\leq 0$, we define it as $+\infty$ in that case to ensure that all iterates remain positive definite. Such a strategy could, however, be detrimental to a solver like {\sc fista}, and to avoid this we selected the starting point so that the condition $\det (P ) \leq 0$  did not occur.

We employ three data sets: the well-known {\tt Estrogen} and {\tt Leukemia} test sets  \cite{li2010inexact}, and the problem given in Olsen et al. \cite{olsen2012newton}, which we call {\tt OONR}. The characteristics of the data sets are given in Table~1, where  nnz($P^\ast_\mu$) denotes the number of nonzeros in the solution. 
 \begin{center}
 Table 1.\\
\begin{tabular}{|r | c c c |}
\hline
Data set & number of features & $\mu$ & nnz($P^\ast_\mu$)\\
\hline
Estrogen & 692 & 0.5 & 10,614 (2.22\%) \\
Leukemia & 1255 & 0.5 & 34,781 (2.21\%) \\
OONR & 500 & 0.5 & 1856 (0.74\%) \\
\hline
\end{tabular}
\end{center}

\bigskip
The performance of the algorithms  on these three test problems is given in Tables~2, 3 and 4. We note that {\sc  fista} does not perform inner iterations since it is applied directly to the original problem \eqref{l1prob}, and that {\sc pnopt-fista} does not compute Hessian vector products because the matrix $H_k$ in the model \eqref{lassop} is defined by quasi-Newton updating. Each inner iteration of  {\sc sl-obm-qn} performs a Hessian-vector multiplication to compute the subproblem objective, and a multiplication of the inverse Hessian times a vector to compute the unconstrained minimizer on the active orthant face --- we report these as two Hessian-vector products in  Tables~2, 3 and 4.

\bigskip


\small
\begin{center}
{\large Table 2.} {ESTROGEN}; $\mu = 0.5$, optimality tolerance = $10^{-5}$ \\ 
\begin{tabular}{|r | c c c c  c| }
\hline 
solver       &  FISTA &  SQA   & PNOPT & SQA &   SQA \\ \hline\hline
inner solver &       & {\sc fista} & {\sc fista} & {\sc obm-qn}   &{\sc obm-cg} \\
\hline
  outer iterations & 808 & 9 & 43 & 44 & 8 \\
  inner iterations & - & 183 & 2134$^\ast$ & 64 & 93 \\
  function/gradient evals & 1751 & 10 & 44 & 45 & 10 \\
  Hessian-vect mults & - & 417 & - & 2 & 213 \\
  time (s) & 208.87 & 51.54 & 355.15 & 38.74 & 26.95 \\
\hline
\multicolumn{6}{l}{$^{\ast}$ For PNOPT we report the number of prox. evaluations}\\
\end{tabular}
\end{center}
\normalsize

\bigskip

\small
\begin{center}
{\large Table 3.} {LEUKEMIA}; $\mu = 0.5$, optimality tolerance = $10^{-5}$ \\ 
\begin{tabular}{|r | c c c c  c |}
\hline 
solver       & FISTA &  SQA   & PNOPT$^{\ast}$ & SQA$^{\ast}$ &  SQA \\ \hline\hline
inner solver &       & {\sc fista} & {\sc fista} & {\sc obm-qn}  &  {\sc obm-cg} \\
\hline
 outer iterations & 838 & 8 & $>488^{\ast\ast}$ & 101 & 8 \\
 inner iterations & - & 187 & - & 196 & 101 \\
 function/gradient evals & 1803 & 9 & - & 103 & 9 \\
 Hessian-vect mults & - & 420 & - & 4 & 239 \\
 time (s) & 1048.77 & 239.23 & - & 171.41 & 140.33 \\
\hline
\multicolumn{6}{l}{$^{\ast}$ out of memory for memory size = 50, we decrease memory size to 5}\\
\multicolumn{6}{l}{$^{\ast\ast}$ exit with message: ``Relative change in function value below ftol''}\\
\multicolumn{6}{l}{$^{\ast\ast}$ optimality error is below $1e-4$ after iteration 73, it is $2.3136e-05$ at termination}\\
\end{tabular}
\end{center}
\normalsize

\bigskip

\small
\begin{center}
{\large Table 4.} {OONR}; $\mu = 0.5$, optimality tolerance = $10^{-5}$ \\ 
\begin{tabular}{|r | c c c c  c |}
\hline
solver       & FISTA &  SQA       &  PNOPT      & SQA  &  SQA \\ \hline\hline
inner solver &       & {\sc fista} &  {\sc fista}  & {\sc obm-qn} &  {\sc obm-cg} \\
\hline
  outer iterations & 212 & 10 & 39 & 37 &  7 \\
  inner iterations & $-$ & 80 & 761 & 37 &  60 \\
  function/gradient evals & 461 & 11 & 41 & 44 &  9 \\
  Hessian-vect mults & $-$ & 193 & - & 2 &  125 \\
  time (s) & 23.53 & 10.14 & 70.37 & 12.73 &  7.09 \\
\hline
\end{tabular}
\end{center}
\normalsize

\bigskip
We now comment on the results given in Tables~2-4. For the inverse covariance selection problem \eqref{cov}, Hessian-vector products are not as expensive as for other problems (c.f. Tables 5-6) --- in fact, these products are not much costlier  than computations with the limited memory BFGS matrix. This fact, combined with the effectiveness of the {\sc obm} method, makes  {\sc sqa-obm-cg} the most efficient of the methods tested. {\sc obm} is a good subproblem solver due to its ability to estimate the set of zero variables quickly, so that the subspace step is computed in a small reduced space  (the density of $P_\mu^\ast$ is less than $2.5\%$ for the three test problems.) In addition, the {\sc obm-cg} method can decrease $\| F_q \|$ drastically in a single iteration, often yielding a high quality {\sc sqa} step and thus  a low number of outer iterations.

We note that the quasi-Newton algorithms {\sc sl-obm-qn} and {\sc pnopt} are different methods because of the subproblem solvers they employ. {\sc sl-obm-qn} uses the two-phase {\sc obm} method in which the quasi-Newton step is computed in a subspace, whereas {\sc pnopt} applies the {\sc fista} iteration to subproblem \eqref{quadm} where $H_k$ is a quasi-Newton matrix.  Although the number of outer iterations of both methods is comparable for problems {\tt Estrogen} and {\tt OONR}, there is a large difference in the number of inner iterations due to power of the {\sc obm} approach.

Note that the number of inner {\sc fista} iterations in {\sc sqa-fista} is always smaller than for {\sc fista}. We repeated the experiment with problem {\tt OONR} using looser optimality tolerances (TOL); the total number of {\sc fista} is given in 
Table~5.
\bigskip

\small
\begin{center}
{\large Table 5.} Effect of convergence tolerance TOL; OONR  \\ 
\begin{tabular}{|r  |  c c  c | }
\hline
  TOL        & $10^{-2}$    & $10^{-3}$    & $10^{-4}$   \\ \hline\hline
  FISTA (\# of outer iterations) & 30 & 74 & 136  \\
  SQA-FISTA (\# of inner iterations) &  47 & 56 &  69 \\
  \hline
\end{tabular}
\end{center}
\normalsize
These results are typical for the covariance selection problems, where the {\sc sqa-fista} is clearly more efficient than {\sc fista}; we will see that this is not the case for the problems considered next.

\subsection{Logistic Regression Problems} 
In our second set of experiments the function $f$ in \eqref{l1prob} is given by a logistic function.  Given $N$ data pairs $(z_i,y_i)$, with $z_i\in\Re^n, \, y_i\in\{-1,1\},\,  i=1,\dots,N$, the optimization problem is given by
\[ 
\min_x \ \ \frac{1}{N}\sum_{i=1}^{N}\log(1+\exp(-y_ix^Tz_i)) + \mu \|x\|_1.
\]
We employed the data given in Table~6, which was downloaded from the SVMLib repository. The values of the regularization parameter $\mu$ were taken from Lee et al. \cite{lee2012proximal}.

\bigskip
\begin{center}
{\large Table 6. Test problems for logistic regression tests}
\begin{tabular}{|r | c c c c |}
\hline
Data set  & $N$ & number of features & $\mu$ & nnz($x^\ast_\mu$)\\
\hline
Gisette (scaled) & 6,000 & 5,000 & 6.67e-04 & 482 (9.64\%) \\
RCV1 (binary) & 20,242 & 47,236 & 3.31e-04 & 140 (0.30\%)\\
\hline
\end{tabular}
\end{center}

\bigskip

\small
\begin{center}
{\large Table 7.} {GISETTE}; $\mu = 6.67e-04$, optimality tolerance = $10^{-5}$ \\ 
\begin{tabular}{|r | c c c c  c |}
\hline
solver       & FISTA & SQA & PNOPT         & SQA    &  SQA \\ \hline\hline
inner solver &      & {\sc fista} & {\sc fista}  & {\sc obm-qn}  &  {\sc obm-cg} \\
\hline
  outer iterations & 1023 & 11 & 237 & 253 &  10\\
  inner iterations & $-$ & 1744 & 25260$^\ast$ & 1075 &  770\\
  function/gradient evals & 2200 & 12 & 240 & 254 &  11\\
  Hessian-vector mults & $-$ & 3761 & $-$ & 2 &  3321\\
  time & 185.55 & 311.28 & 108.84 & 38.47 &  273.10 \\
\hline
\multicolumn{6}{l}{$^\ast$ For PNOPT we report the number of prox. evaluations.}
\end{tabular}
\end{center}
\normalsize

\bigskip

\small
\begin{center}
{\large Table 8.} {RCV1}; $\mu = 3.31e-04$, optimality tolerance = $10^{-5}$ \\ 
\begin{tabular}{|r | c c c c  c | }
\hline

solver & FISTA & SQA & PNOPT & SQA &  SQA \\ \hline\hline
inner solver & & {\sc fista} & {\sc fista} &  {\sc obm-qn}  &  {\sc obm-cg} \\
\hline
  outer iterations & 90 & 7 & 19 & 18 & 6 \\
  inner iterations & $-$ & 366 & 1148 & 27 & 54 \\
  function/gradient evals & 184 & 8 & 20 & 19 & 7 \\
  Hessian-vector mults & $-$ & 738 & $-$ & 2 & 120 \\
  time (s) & 1.95 & 7.52 & 11.23 & 0.92 &  1.33 \\
\hline
\end{tabular}
\end{center}

\bigskip\normalsize
For the logistic regression problems,  Hessian-vector products are expensive, particularly for {\tt gisette}, where the data set is dense. As a result, the {\sc obm} variant that employs quasi-Newton approximations, namely {\sc sqa-obm-qn}, performs best (even though {\sc sqa-obm-cg} requires a smaller number of outer iterations). Note that {\sc sqa-fista} is not efficient; in fact it requires a much larger number of inner iterations than the total number of iterations in {\sc fista}. In Table~9 we observe the effect of the optimality tolerance, on these two methods, using problem {\tt gisette}.
\bigskip

\small
\begin{center}
{\large Table 9.} Effect of convergence tolerance TOL ; Gisette \\ 
\begin{tabular}{|r  |  c c  c | }
\hline
  TOL        & $10^{-4}$    & $10^{-5}$    & $10^{-6}$   \\ \hline\hline
  FISTA (\# of outer iterations) & 605 & 1023 & 2555  \\
  SQA-FISTA (\# of inner iterations) &  1249 & 1744 &  2002 \\
  \hline
\end{tabular}
\end{center}
\normalsize
\bigskip

We observe from Table~9, that {\sc fista} requires a smaller number of iterations; it is only for a very high accuracy of
$10^{-6}$ that {\sc sqa-fista} becomes competitive. This is in stark contrast with Table~5.

In summary, for the logistic regression problems the advantage of the {\sc sqa} method is less pronounced than for the inverse covariance estimation problems, and is achieved only through the appropriate choice of model Hessian $H_k$ (quasi-Newton) and the appropriate choice of inner solver (active set {\sc obm} method).

\subsection{Description of the orthant based method (OBM) } \label{obm}
We conclude this section by describing the orthant-based method  used in our experiments to solve the subproblem \eqref{lassop}. We let $t$ denote the iteration counter of the {\sc obm} method, and let $z_t$ denote its iterates.

Given an iterate $z_t$, the method defines an orthant face $\Omega_t$ of $\R^n$ by
 \begin{equation}    \label{ao}
    \Omega_t = \mbox{\textbf{cl}}(\{ d \in \R^n: \sgn(d_i) = \sgn([\omega_t]_i), \ i=1,...,n \}),   
 \end{equation}
 with
 \begin{equation}
 [\omega_t]_i = \begin{cases} \sgn([z_t]_i) \quad \ \ \mbox{if }  [z_t]_i\neq 0 \\ \sgn(-[v_t]_i) \quad \mbox{if } [z_t]_i= 0 , \end{cases}  
\end{equation}
 where $v_t$ is the minimum norm subgradient of $q_k$ computed at $z_t$, i.e.,
    \begin{equation}    \label{mnsub}
      [v_t]_i = \begin{cases} [\nabla q_k(z_t)]_i + \mu \quad &\mbox{ if } [z_t]_i> 0 \quad \mbox{or }  \mbox{ }( [z_t]_i = 0 
                           \mbox{ } \land \mbox{ } \nabla q_k(z_t)]_i + \mu < 0 ) \\
                          [\nabla q_k(z_t)]_i - \mu \quad &\mbox{ if } [z_t]_i< 0 \quad \mbox{or } \mbox{ }  ([z_t]_i= 0 
                          \mbox{ } \land \mbox{ } \nabla q_k(z_t)]_i - \mu > 0 ) \\
                           \qquad \quad 0 \quad   \quad &\mbox{ if }  [z_t]_i= 0 
                           \quad \mbox{and }  \mbox{ } 0 \in [\nabla q_k(z_t)]_i - \mu,\nabla q_k(z_t)]_i + \mu] . \\ \end{cases}
 \end{equation}
 Defining $\Omega_t$ in this manner was proposed, among others, by Andrew and Gao \cite{andrew2007scalable}.
 In the relative interior of $\Omega_t$, the model function $q_k$ is differentiable. The active set in the orthant-based method, defined as $A^k = \{i : \omega^k_i =0\}$, determines the variables that are kept at zero, while the rest of the variables are chosen to minimize a (smooth) quadratic model. Specifically, the search direction $d_t$  of the algorithm is given by $d_t=\hat z -z_t$, where $\hat z$ is a solution of 
 \begin{align}    
      \min_{z \in \R^{n}} & \ \ \psi(z) =   q_k(z_t) + (z-z_t)^Tv^k 
                                          + \fraction{1}{2} (z-z_t)^T  H_k (z-z_t) \nonumber \\
      \mbox{s.t.} & \ \ z_i=[z_t]_i, \ i \in A^k . \label{aomodel}
\end{align}  
Note that $\psi(z) = f(x_k) + (g(x_k)+\omega_t\mu)^T(z-x_k) + \frac{1}{2}(z-x_k)^TH_k(z-x_k)$.

In the {\sc obm-cg} variant, we set $H_k = \nabla^2  f(x_k)$, and perform an  approximate minimization of this problem using the projected conjugate gradient iteration \cite{mybook}. In the {\sc obm-qn} version, $H_k$ is a limited memory BFGS matrix and $\hat z$ is the exact solution of \eqref{aomodel}. This requires computation of the inverse reduced Hessian $R_k=(Z_k^T \nabla^2 H_k Z_k)^{-1}$, where $Z_k$ is a basis for the space defined by \eqref{aomodel}. The matrix $R_k$ can be updated using  the compact representations of quasi-Newton matrices \cite{ByrdNoceSch94}.
After the direction $d_t= \hat z-z_t$ has been computed, the {\sc obm} method performs a line search along $d_t$, projecting the iterate back onto the orthant face $\Omega_t$, until a sufficient reduction in the function $q_k$ has been obtained.  Although this algorithm performed reliably in our tests, its convergence has not been proved (to the best of our knowledge) because the orthant face identification procedure \eqref{ao}-\eqref{aomodel} can lead to arbitrarily small steps.

Our {\sc obm-qn} algorithm differs from the {\sc owl} method in two respects: it does not realign the direction $z-z_t$ so that the sign of its components match those of $v_t$, and it performs the minimization of the model exactly, while the {\sc owl} method computes only an approximate solution -- defined by computing the reduced inverse Hessian $Z_k^T \nabla^2 H_k^{-1} Z_k$,
instead of the inverse of the reduced Hessian $R_k$.

\section{Final Remarks} \label{finalr}
One of the key ingredients in making the successive quadratic approximation (or proximal Newton) method practical for problem \eqref{l1prob} is the ability to terminate the inner iteration as soon as a step of sufficiently good quality is computed. In this paper, we have proposed such an inexactness criterion; it employs an optimality measure that is tailored to the structure of the problem. We have shown that the resulting algorithm is globally convergent, that its rate of convergence can be controlled through an inexactness parameter, and that the inexact method will naturally accept unit step lengths in a neighborhood of the solution. We have also argued that our inexactness criterion is preferable to the one proposed by Lee et al. \cite{lee2012proximal}.

The method presented in this paper can use any algorithm for the inner minimization of the subproblem \eqref{quadm}. In particular, all the results are applicable to the case  when this inner minimization is performed using a coordinate descent algorithm \cite{hsieh2011sparse, tanscheinberg}. In our numerical tests we employed {\sc fista} and an orthant-based method as the inner solvers, and found the latter method to be particularly effective. The efficacy of the successive quadratic approximation approach depends of the choice of matrix $H_k$ in \eqref{quadm}, which is problem dependent:  when Hessian-vector products are expensive to compute, then a quasi-Newton approximation is most efficient; otherwise defining $H_k$ as the exact Hessian and implementing a Newton-CG iteration is likely to give the best results.

\bigskip\noindent{\it Acknowledgement.} The authors thank Jong-Shi Pang for his insights and advice throughout the years. The theory presented by Facchinei and Pang \cite{facchinei2003finite} in the context of variational inequalities was used in our analysis,  showing the power and generality  that masterful book.

\small
\bibliographystyle{plain}
\bibliography{../../References/references}
\end{document}